\documentclass[11pt]{article}
\usepackage{graphicx}
\usepackage{amsmath}
\DeclareMathOperator{\arctanh}{arctanh}
\usepackage{amssymb}
\usepackage{amsthm}
\usepackage[titletoc]{appendix}
\usepackage{url}

\newcommand\al\alpha
\newcommand\be\beta
\newcommand\de\delta
\newcommand\ep\varepsilon
\newcommand\tha\theta
\newcommand\ka\kappa
\newcommand\la\lambda
\newcommand\om\omega

\newcommand\iy\infty
\newcommand\pa\partial

\renewcommand\Re{\operatorname{Re}}
\numberwithin{equation}{section}

\newtheorem{theorem}{Theorem}

\newtheorem{Remark}[theorem]{Remark}

\begin{document}
	
\title{Some integrals which are not in table 129 of Bierens de Haan.}
\author{Enno Diekema \footnote{email address: e.diekema@gmail.com}}
\maketitle
	
\begin{abstract}
\noindent
In the standard work of Bierens de Haan about integrals we look at table 129. This table lists a number of integrals of a certain kind. In this paper the table is expanded with a number of similar integrals. These are determined by a number of different methods. Some of these methods are not treated in the freshman standard works on integral calculus. As a bonus, we treat an integral using a complete different method which is not so well-known.
\end{abstract}
	
\section{Introduction}
\setlength{\parindent}{0cm}
\noindent
	
During my investigations on the function $\zeta(3)$ I came across a number of integrals that resembles some integrals that appear in the standard tables of Bierens de Haan \cite{1}. This concerns table 129 in \cite{1}. A copy of this table is in the Appendix. A number of integrals in this table are connected with table 97 of \cite{1}. Most of these integrals have been adopted in the work of Gradshteyn and Ryzhik \cite{2}. Moll \cite{7} proves a number of these integrals with different methods. Some of these methods we use in this paper.

However, we must note that there is a large number of errors in the tables of Bierens de Haan. See \cite{3}. This means that every integral taken from \cite{1} must be checked to see if it is correct.

In this paper we distinguish four basic integrals. These are respectively
\[
\qquad \int_0^1\dfrac{1}{a^2+\ln(x)^2}\dfrac{1}{1+x}dx 
\qquad \int_0^1\dfrac{1}{a^2+\ln(x)^2}\dfrac{1}{(1+x)^2}dx 
\]
\[
\qquad \int_0^1\dfrac{\ln(x)}{a^2+\ln(x)^2}\dfrac{1}{1+x}dx 
\qquad \int_0^1\dfrac{\ln(x)}{a^2+\ln(x)^2}\dfrac{1}{(1+x)^2}dx 
\]
The integrals 1..11 belong to the group of the first integral. The integrals 12 and 13 belong to the group of the second integral. The integrals 14..17 belong to the group of the third integral and the integral 18 belongs to the group of the fourth integral. 

By differentiating the integrals in the four groups with respect to the variable $a$ complete classes of integrals arise. Special cases occur when $a=\pi^2$. The last integral 19 does not belong to the group of the other integrals. We add this integral because we use two methods to get the answer. The first method use the differentiation with respect to the variable $a$. The second method uses a totally different method which is not often used.

\

For the computation of the integrals we use a number of methods. The first method uses power series, the second method is due to Schl\"omilch, the third method is due to Saalsch\"utz. These last two methods are not in the common textbooks. Note that different methods gives mostly different results. 

\

In section 2 there is an overview of the results. In section 3 we give all the special functions with some properties we need to derive the integrals. In section 4 we prove all the integrals from section 3. In section 5 we have a remarkable overview of a number of integrals. In the Appendix we show the integrals from Table 129 from Bierens de Haan \cite{1}.

\section{Main results}
The main results of this paper are the following integrals
\begin{align*}
&1.\quad\int_0^1\dfrac{1}{a^2+\ln(x)^2}\dfrac{1}{1+x}dx=
\dfrac{1}{a}\sum_{k=1}^\infty(-1)^k\left(\left(Si(a\, k)-\dfrac{\pi}{2}\right)\cos(a\, k)
-Ci(a\, k)\sin(a\, k)\right) \\
&2.\quad\int_0^1\dfrac{1}{\pi^2+\ln(x)^2}\dfrac{1}{1+x}dx=\dfrac{1}{\pi}\sum_{k=1}^\infty \text{si}\big(k\pi\big) \\
&3.\quad\int_0^1\dfrac{1}{a^2+\ln(x)^2}\dfrac{1}{1+x}dx=\dfrac{\pi}{4a}-2\sum_{k=0}^\infty\dfrac{\ln(a)-\ln(\pi)-\ln(2k+1)}{a^2-(2k+1)^2\pi^2} \\
&4.\quad\int_1^\infty\dfrac{1}{a^2+\ln(x)^2}\dfrac{1}{1+x}dx=\dfrac{\pi}{4a}+2\sum_{k=0}^\infty\dfrac{\ln(a)-\ln(\pi)-\ln(2k+1)}{a^2-(2k+1)^2\pi^2} \\
&5.\quad\int_0^\infty\dfrac{1}{a^2+\ln(x)^2}\dfrac{1}{1+x}dx=\dfrac{\pi}{2a} \\
&6.\quad\int_0^1\dfrac{1}{\pi^2+\ln(x)^2}\dfrac{1}{1+x}dx=
\dfrac{1}{4}-\dfrac{1}{\pi^2}-\dfrac{1}{2\pi^2}\sum_{k=1}^\infty\dfrac{\ln(2k+1)}{k(k+1)} \\
&7.\quad\int_0^1\dfrac{1}{\pi^2+\ln(x)^2}\dfrac{1}{1+x}dx=
\dfrac{1}{4}-\dfrac{1}{\pi^2}-\dfrac{1}{\pi^2}\ln\left(\dfrac{\pi}{2}\right)+
\dfrac{1}{2}\sum_{k=1}^\infty\dfrac{(-1)^k\pi^{2k-2}}{(2k)!k(2k-1)}B_{2k} \\
&8.\quad\int_0^1\dfrac{1}{\pi^2+\ln(x)^2}\dfrac{1}{1+x}dx=
\dfrac{1}{4}-\dfrac{1}{\pi^2}-\dfrac{2}{\pi^2}\sum_{k=1}^\infty\dfrac{1}{(2k-1)}\dfrac{1}{2^{2k}}
\zeta(2k) \\
&9.\quad\int_0^1\dfrac{1}{4\pi^2+\ln(x)^2}\dfrac{1}{1+x}dx=
\dfrac{1}{8}-\dfrac{2}{\pi^2}\sum_{k=1}^\infty\dfrac{\ln(2k+1)}{(2k-1)(2k+3)} \\
&10.\quad\int_0^1\dfrac{1}{\left(a^2+\ln(x)^2\right)^2}\dfrac{1}{1+x}dx=
\dfrac{\pi}{8a^3}-\dfrac{1}{4a^3}\tan\left(\dfrac{a}{2}\right)+
2\sum_{k=0}^\infty\dfrac{\ln\big((2k+1)\pi\big)-\ln(a)}{\big((1+2k)^2\pi^2-a^2\big)^2} \\
&11.\quad\int_0^1\dfrac{1}{\left(\pi^2+\ln(x)^2\right)^2}\dfrac{1}{1+x}dx=
\dfrac{1}{8\pi^2}-\dfrac{3}{4\pi^4}+\dfrac{1}{8\pi^4}\sum_{k=1}^\infty\dfrac{\ln(1+2k)}{k^2(1+k)^2} \\
&12.\quad\int_0^1\dfrac{1}{a^2+\ln(x)^2}\dfrac{1}{(1+x)^2}dx=
-\dfrac{1}{4a\pi}\psi^{(1)}\left(\dfrac{a}{2\pi}\right)+
\dfrac{1}{a\pi}\psi^{(1)}\left(\dfrac{a}{\pi}\right) \\
&13.\quad\int_0^1\dfrac{1}{\pi^2+\ln(x)^2}\dfrac{1}{(1+x)^2}dx=
\int_1^\infty\dfrac{1}{\pi^2+\ln(x)^2}\dfrac{1}{(1+x)^2}dx=
\dfrac{1}{\pi}\sum_{k=1}^\infty k\, \text{si}(k\pi)=\dfrac{1}{24} \\
&14.\quad\int_0^1\dfrac{\ln(x)}{a^2+\ln(x)^2}\dfrac{1}{1+x}dx=
\dfrac{1}{2}\ln\left(\dfrac{2a}{\pi}\right)+
\dfrac{1}{2}\psi\left(\dfrac{a}{2\pi}\right)-\psi\left(\dfrac{a}{\pi}\right) 
\end{align*}
\begin{align*}
&15.\quad\int_0^1\dfrac{\ln(x)}{\pi^2+\ln(x)^2}\dfrac{1}{1+x}dx=
-\sum_{k=1}^\infty\text{Ci}(k\, \pi)=\dfrac{\gamma}{2}-\dfrac{1}{2}\ln(2) \\
&16.\quad\int_0^1\dfrac{\ln(x)}{\left(a^2+\ln(x)^2\right)^2}\dfrac{1}{1+x}dx=
-\dfrac{1}{4a^2}-\dfrac{1}{8a\pi}\psi^{(1)}\left(\dfrac{a}{2\pi}\right)+
\dfrac{1}{2a\pi}\psi^{(1)}\left(\dfrac{a}{\pi}\right) \\
&17.\quad\int_0^1\dfrac{\ln(x)}{\left(\pi^2+\ln(x)^2\right)^2}\dfrac{1}{1+x}dx=
\dfrac{1}{48}-\dfrac{1}{4\pi^2} \\
&18.\quad\int_0^1\dfrac{\ln(x)}{\pi^2+\ln(x)^2}\dfrac{1}{(1+x)^2}dx=
-\sum_{k=1}^\infty k\ \text{Ci}(k\, \pi) \\
&19.\quad\int_0^1\dfrac{1}{\big(\pi^2+\ln(x)^2\big)^2}\dfrac{1}{(x+1)^2}dx=
\int_1^\infty\dfrac{1}{\big(\pi^2+\ln(x)^2\big)^2}\dfrac{1}{(x+1)^2}dx=
\dfrac{\zeta(3)+\zeta(2)}{8\pi^4} 
\end{align*}

\section{Overview of some formulas of the used special functions}
\subsection{Gamma function}
The Gamma function is defined as
\[
\Gamma(z)=\int_0^\infty t^{z-1}e^{-t}dt \qquad \Re{(z)}>0
\]
For $z=n$ we get $\Gamma(n)=(n-1)!$
\subsection{Polygamma function}
The polygamma function of order $n$ is defined as the $(n + 1)$th derivative of the logarithm of the gamma function
\[
\psi^{(n)}(z)=\psi(n,z)=\dfrac{d^n}{dz^n}\psi(z)=
\]
For $n=1$ we get
\[
\psi^{(1)}(z)=\psi(z)=\dfrac{d}{dz}\ln\Gamma(z)=\dfrac{\Gamma'(z)}{\Gamma(z)}
\]
$\psi(z)$ is called the digamma function. 

A well-known property \cite[8.365.8]{2} is
\[
\psi(1-z)-\psi(z)=\pi\cot(\pi\, z) 
\]
\subsection{Dilogarithm function}
The polylogarithm function is defined as
\[
Li_2(z)=-\int_0^z\dfrac{\ln(1-t)}{t}dt
\]

\subsection{Euler-Mascheroni constant}
The Euler-Mascheroni constant $\gamma$ and is given by
\[
\gamma=\lim_{n \rightarrow \infty}\left(\sum_{k=1}^n\dfrac{1}{k}-\ln(n)\right)
=-\int_0^\infty\dfrac{\ln(x)}{e^x}dx
\]
\subsection{Cosine integral function}
The cosine integral function is given by
\[
\text{Ci(z)}=-\int_z^\infty\dfrac{\cos(t)}{t}dt=\gamma+\ln(x)+\int_0^x \dfrac{\cos(t)-1}{t}dt
\]
\subsection{Sine integral functions}
The sine integral function is given by
\[
\text{Si(z)}=\int_0^z\dfrac{\sin(t)}{t}dt
\]
The next form is also often used
\[
\text{si}(z)=-\int_z^\infty\dfrac{\sin(t)}{t}dt
\]
There follows
\[
\text{si}(z)=\text{Si}(z)-\dfrac{\pi}{2}
\]
\subsection{Exponential integral}
The exponential integral is given by
\[
\text{Ei}(x)=-\int_{-x}^\infty\dfrac{e^{-t}}{t}dt=\int_{-\infty}^x\dfrac{e^t}{t}dt
\]
We need the following property
\begin{equation}
\dfrac{i}{2}\Big[\text{Ei}(-ix)-\text{Ei}(ix)\Big]=\text{Si}(x)+\dfrac{\pi}{2}
\label{s2.6}
\end{equation}
\subsection{Bernoulli numbers}
The Bernoulli numbers are given by
\begin{equation}
B_m=\sum_{k=0}^m\dfrac{1}{k+1}\sum_{j=0}^k\binom{k}{j}(-1)^j(j+1)^m \qquad m \geq 1
\label{s2.7}
\end{equation}
\subsection{Riemann zeta function}
The Riemann zeta function is given by
\[
\zeta(s)=\sum_{n=1}^\infty\dfrac{1}{n^s}=
\dfrac{1}{\Gamma(s)}\int_0^\infty\dfrac{x^{s-1}}{e^x-1}dx
\]
If $s$ is an integer $2k$ then there is a connection between the Bernoulli number and the zeta number
\[
B_{2k}=\dfrac{(-1)^{k+1}2(2k)!}{(2\pi)^{2k}}\zeta(2k)
\]

\section{Proof of the integrals}

\subsection{Proof of integral 1}
In this section we treat integral $1$
\[
I=\int_0^1\dfrac{1}{a^2+\ln(x)^2}\dfrac{1}{1+x}dx
\]
Because $x$ lies between $0$ and $1$ the second factor in the integral can be written as a geometric series
\[
\dfrac{1}{1+x}=\sum_{k=0}^\infty(-x)^k 
\]
The integral becomes
\[
I=\int_0^1\dfrac{1}{a^2+\ln(x)^2}\dfrac{1}{1+x}dx=\int_0^1\dfrac{1}{a^2+\ln(x)^2}\sum_0^\infty(-x)^kdx
\]
Because both factors are absolute convergent we may interchange the summation and the integral
\begin{equation}
I=\int_0^1\dfrac{1}{a^2+\ln(x)^2}\dfrac{1}{1+x}dx=
\sum_{k=0}^\infty\int_0^1\dfrac{(-x)^k}{a^2+\ln(x)^2}dx
\label{3.1}
\end{equation}
The integral is known and we get
\[
\int_0^1\dfrac{(-x)^k}{a^2+\ln(x)^2}dx=
\dfrac{1}{a}(-1)^k\left(\text{Ci(a(1+k))}\sin(a(1+k))-\left(\text{Si(a(1+k) )}
+\dfrac{\pi}{2}\right)\cos(a(1+k))\right)
\]
Substitution in \eqref{3.1} gives
\[
I=\dfrac{1}{a}\sum_{k=0}^\infty(-1)^k\left(\text{Ci(a(1+k))}\sin(a(1+k))-\left(\text{Si(a(1+k) )}
+\dfrac{\pi}{2}\right)\cos(a(1+k))\right)
\]
Starting the summation with $k=1$ we get
\[
\int_0^1\dfrac{1}{a^2+\ln(x)^2}\dfrac{1}{1+x}dx=
\dfrac{1}{a}\sum_{k=1}^\infty(-1)^k\left(\left(Si(a\, k)-\dfrac{\pi}{2}\right)\cos(a\, k)
-Ci(a\, k)\sin(a\, k)\right)
\]
and this completes the proof.

\subsection{Proof of integral 2}
In this section we treat integral $2$ 
\[
I=\int_0^1\dfrac{1}{\pi^2+\ln(x)^2}\dfrac{1}{1+x}dx
\]
We use two methods. The first method uses the result of the previous integral. The second method uses a Feynmann trick.

\subsubsection{Direct proof}
For proving this integral we set $a=\pi$ in integral 1.
\[
I=\int_0^1\dfrac{1}{\pi^2+\ln(x)^2}\dfrac{1}{1+x}dx=
\dfrac{1}{\pi}\sum_{k=1}^\infty(-1)^k\left(\left(\text{Si}(\pi\, k)
-\dfrac{\pi}{2}\right)\cos(\pi\, k)-\text{Ci}(\pi\, k)\sin(\pi\, k)\right)
\]
Simplification gives
\begin{equation}
I=\int_0^1\dfrac{1}{\pi^2+\ln(x)^2}\dfrac{1}{1+x}dx=
\dfrac{1}{\pi}\sum_{k=1}^\infty\left(Si(\pi\, k)-\dfrac{\pi}{2}\right)
\label{3.2}
\end{equation}
Using the definition of the function $\text{si}(z)$ we get
\[
I=\int_0^1\dfrac{1}{\pi^2+\ln(x)^2}\dfrac{1}{1+x}dx=
\dfrac{1}{\pi}\sum_{k=1}^\infty \text{si}\big(k\pi\big)
\]
which completes the proof.

\subsubsection{Using a Feynman trick}
Feynman uses his well-known trick to calculate integrals. To this end, a parameter is added somewhere in the integral. The integral can be calculated by applying a number of operations to the parameter.

We start with applying a parameter $b$ in the second fraction
\[
\int_0^1\dfrac{1}{\pi^2+\ln(x)^2}\dfrac{1}{1+x}dx \rightarrow 
\int_0^1\dfrac{1}{\pi^2+\ln(x)^2}\dfrac{1}{1+b\, x}dx
\]
After integration and differentiation to the parameter $b$ we get after interchanging the integral and the differential operator (which can be proved that these operations are permitted)
\begin{align*}
\int_0^1\dfrac{1}{\pi^2+\ln(x)^2}\dfrac{1}{1+b x}dx
&=\int_0^1\dfrac{1}{\pi^2+\ln(x)^2}\left(\dfrac{d}{db}\int\dfrac{1}{1+b\,  x}db\right)dx \\
&=\dfrac{d}{db}\int_0^1\dfrac{1}{\pi^2+\ln(x)^2}\dfrac{\ln(1+b\, x)}{x}dx
\end{align*}
The last fraction can be written as a power series
\begin{align*}
\int_0^1\dfrac{1}{\pi^2+\ln(x)^2}\dfrac{1}{1+b x}dx
&=\dfrac{d}{db}\sum_{k=1}^\infty(-1)^{k-1}\dfrac{b^k}{k}\int_0^1\dfrac{1}{\pi^2+\ln(x)^2}
x^{k-1}dx \nonumber \\
&=\sum_{k=1}^\infty(-1)^{k-1}b^{k-1}\int_0^1\dfrac{1}{\pi^2+\ln(x)^2}
x^{k-1}dx
\end{align*}
For the integral of $x$ we get
\[
\int_0^1\dfrac{1}{\pi^2+\ln(x)^2}x^{k-1}dx=
\dfrac{(-1)^k}{2\pi}
\big(2\pi-i\left(\text{Ei}(-ik\pi)-\text{Ei}(ik\pi)\right)\big)
\]
Using \eqref{s2.6}, setting $b=1$ and gathering all the previous results gives at last
\[
\int_0^1\dfrac{1}{\pi^2+\ln(x)^2}\dfrac{1}{1+x}dx=
\dfrac{1}{\pi}\sum_{k=1}^\infty\left(\text{Si}(k\pi)-\dfrac{\pi}{2}\right)
\]
This is the same as equation \eqref{3.2} and the proof is complete.

\subsection{Proof of integral 3}
In this section we treat again integral $3$
\[
\int_0^1\dfrac{1}{a^2+\ln(x)^2}\dfrac{1}{1+x}dx
\]
with a formula of Saalsch\"utz. Transformation of the integral with $x=\exp(-t)$ gives
\[
I=\int_0^1\dfrac{1}{a^2+\ln(x)^2}\dfrac{1}{1+x}dx=
\int_0^\infty\dfrac{1}{a^2+t^2}\dfrac{\exp(-t/2)}{\exp(t/2)+\exp(-t/2)}dt
\]
Rewriting the last fraction gives
\[
\dfrac{\exp(-t/2)}{\exp(t/2)+\exp(-t/2)}=
\dfrac{1}{2}-\dfrac{1}{2}\left(\dfrac{\exp(t/2)-\exp(-t/2)}{\exp(t/2)+\exp(-t/2)}\right)
\]
For the fraction on the righthand side we use the last formula on page 111 of the book of Saalsch\"utz \cite{9}. This formula gives a power series.
\begin{equation}
\dfrac{\exp(x)-\exp(-x)}{\exp(x)+\exp(-x)}=8x\sum_{k=0}^\infty\dfrac{1}{(2k+1)^2\pi^2+4x^2}
\label{3.4.00}
\end{equation}
Application gives
\begin{align*}
I=\int_0^1\dfrac{1}{a^2+\ln(x)^2}\dfrac{1}{1+x}dx
&=\int_0^\infty\dfrac{1}{a^2+t^2}\left(\dfrac{1}{2}-2t\sum_{k=0}^\infty\dfrac{1}{(2k+1)^2\pi^2+t^2}\right)dt \\
&=\dfrac{1}{2}\int_0^\infty\dfrac{1}{a^2+t^2}dt-2\int_0^\infty\dfrac{t}{a^2+t^2}\sum_{k=0}^\infty
\dfrac{1}{(2k+1)^2\pi^2+t^2}dt \\
&=\dfrac{1}{2}\int_0^\infty\dfrac{1}{a^2+t^2}dt-2\sum_{k=0}^\infty\int_0^\infty\dfrac{t}{a^2+t^2}
\dfrac{1}{(2k+1)^2\pi^2+t^2}dt
\end{align*}
After integration we get at last
\[
\int_0^1\dfrac{1}{a^2+\ln(x)^2}\dfrac{1}{1+x}dx=\dfrac{\pi}{4a}-2\sum_{k=0}^\infty\dfrac{\ln(a)-\ln(\pi)-\ln(2k+1)}{a^2-(2k+1)^2\pi^2}
\]
which completes the proof.

\subsection{Proof of integral 4}
In this section we treat integral $4$
\[
\int_1^\infty\dfrac{1}{a^2+\ln(x)^2}\dfrac{1}{1+x}dx
\]
Transformation of the integral with $x \rightarrow \dfrac{1}{x}$ gives
\begin{align*}
\int_1^\infty\dfrac{1}{a^2+\ln(x)^2}\dfrac{1}{1+x}dx
&=\int_0^1\dfrac{1}{a^2+\ln(x)^2}\left(\dfrac{1}{x}-\dfrac{1}{1+x}\right)dx= \\
&=\dfrac{\pi}{2a}-\int_0^1\dfrac{1}{a^2+\ln(x)^2}\dfrac{1}{1+x}dx
\end{align*}
Using integral 4 gives
\begin{align*}
\int_1^\infty\dfrac{1}{a^2+\ln(x)^2}\dfrac{1}{1+x}dx
&=\dfrac{\pi}{2a}-\left(\dfrac{\pi}{4a}-2\sum_{k=0}^\infty\dfrac{\ln(a)-\ln(\pi)-\ln(2k+1)}{a^2-(2k+1)^2\pi^2}\right) \\
&=\dfrac{\pi}{4a}+2\sum_{k=0}^\infty\dfrac{\ln(a)-\ln(\pi)-\ln(2k+1)}{a^2-(2k+1)^2\pi^2}
\end{align*}
This completes the proof.

\subsection{Proof of integral 5}
In this section we treat integral $5$
\[
\int_0^\infty\dfrac{1}{a^2+\ln(x)^2}\dfrac{1}{1+x}dx
\]
Adding integral 3 and integral 4 gives 
\[
\int_0^\infty\dfrac{1}{a^2+\ln(x)^2}\dfrac{1}{1+x}dx=\dfrac{\pi}{2a}
\]
This completes the proof.

\subsection{Proof of integral 6}
In this section we treat integral $6$
\[
\int_0^1\dfrac{1}{\pi^2+\ln(x)^2}\dfrac{1}{1+x}dx
\]
We use the integral 3. For $a=\pi$ this is a special case because the summand  for $k=0$ gives a problem. We have to take the limit. After taking the limit we get
\begin{equation}
\int_0^1\dfrac{1}{\pi^2+\ln(x)^2}\dfrac{1}{1+x}dx=
\dfrac{1}{4}-\dfrac{1}{\pi^2}-\dfrac{1}{2\pi^2}\sum_{k=1}^\infty\dfrac{\ln(2k+1)}{k(k+1)}
\label{3.6}
\end{equation}
This completes the proof.

\subsection{Proof of integral 7}
In this section we treat again integral $6$
\[
\int_0^1\dfrac{1}{\pi^2+\ln(x)^2}\dfrac{1}{1+x}dx
\]
to get another result. 

Using equation \eqref{3.6} we convert the summation with a logarithm into a summation with Bernoulli numbers. The summation with the logarithm is
\[
S=\sum_{k=1}^\infty\dfrac{\ln(2k+1)}{k(k+1)}
\]
In summation problems we can also apply the Feynmann trick. So we write the summation as
\[
S=\sum_{k=1}^\infty\dfrac{\ln(2k+1)}{k(k+1)}=
\sum_{k=1}^\infty\dfrac{\ln(2k+2x)}{k(k+1)}=
\sum_{k=1}^\infty\dfrac{\ln(2)}{k(k+1)}+\sum_{k=1}^\infty\dfrac{\ln(k+x)}{k(k+1)}=
\ln(2)+\sum_{k=1}^\infty\dfrac{\ln(k+x)}{k(k+1)}
\]
At the end of the computation we substitute $x=\dfrac{1}{2}$ into the result. Taking the derivative of the summation with respect to the parameter $x$ gives
\[
T=\sum_{k=1}^\infty\dfrac{\ln(k+x)}{k(k+1)}\qquad\qquad
\dfrac{dT}{dx}=\sum_{k=1}^\infty\dfrac{1}{k(k+1)(k+x)}
\]
We write the factor $\dfrac{1}{k+x}$ as an integral and get
\[
\dfrac{dT}{dx}=\sum_{k=1}^\infty\dfrac{1}{k(k+1)}\int_0^1t^{k+x-1}dt
\]
Interchanging the summation and the integral (it can be shown that this is allowed) gives
\[
\dfrac{dT}{dx}=\int_0^1 t^{x-1}\sum_{k=1}^\infty\dfrac{t^k}{k(k+1)}dt
\]
The summation is well-known and we get
\[
\dfrac{dT}{dx}=\int_0^1 t^{x-1}\big(t+\ln(1-t)-t\ln(1-t)\big)dt
\]
After integration we get
\[
\dfrac{dT}{dx}=-\dfrac{1}{1-x}+\dfrac{1}{(1-x)x}+\dfrac{\gamma}{(1-x)x}+\dfrac{\psi(x)}{(1-x)x}
\]
$\gamma$ is the Euler-Mascheroni constant. Integration gives 
\[
S=\ln(2)+\left[-\dfrac{1}{x}-\gamma\ln(1-x)+(1+\gamma)\ln(x)\right]_{x=1/2}+
\left[\int\dfrac{\psi(x)}{x}dx+\int\dfrac{\psi(x)}{1-x}\right]_{x=1/2}
\]
After some simplification we get
\begin{equation}
S=-2+\left[\int\dfrac{\psi(x)}{x}dx+\int\dfrac{\psi(x)}{1-x}\right]_{x=1/2}
\label{3.4aa}
\end{equation}
With $\dfrac{d}{dx}\psi(x)=\ln\Gamma(x)$ partial integration gives
\[
S=-2+2\ln(\pi)+\left[\int\dfrac{\ln\Gamma(x)}{x^2}dx-\int\dfrac{\ln\Gamma(x)}{(1-x)^2}\right]_{x=1/2}
\]
The first integrand can be written as a power series.
\[
\dfrac{\ln\Gamma(x)}{x^2}=-\dfrac{\gamma}{x}-\dfrac{\ln(x)}{x^2}-
\sum_{k=1}^\infty(-1)^k\dfrac{B_{2k}}{(2k)!~k}2^{2k-2}\pi^{2k}x^{2k-2}+
\sum_{k=1}^\infty\dfrac{1}{(2k+1)!}\psi(2k,1)x^{2k-1}
\]
$\psi(2k,1)$ is the polygamma function. Application gives
\begin{multline*}
\left[\int\dfrac{\ln\Gamma(x)}{x^2}dx\right]_{x=1/2}=
2-2\ln(2)+\gamma\ln(2)-\dfrac{1}{2}\sum_{k=1}^\infty(-1)^k\dfrac{B_{2k}}{(2k)!k(2k-1)}\pi^{2k}+ \\
+\sum_{k=1}^\infty\dfrac{1}{k(2k+1)!}\psi(2k,1)\dfrac{1}{2^{2k+1}}
\end{multline*}
The second integrand can also be written as a power series.
\[
\dfrac{\ln\Gamma(x)}{(1-x)^2}=\dfrac{\gamma}{1-x}-\sum_{k=1}^\infty(-1)^k\dfrac{B_{2k}}{(2k)!~k}2^{2k-2}\pi^{2k}(1-x)^{2k-2}-\sum_{k=1}^\infty\dfrac{1}{(2k+1)!}\psi(2k,1)(1-x)^{2k-1}
\]
Application gives
\[
\left[\int\dfrac{\ln\Gamma(x)}{(1-x)^2}dx\right]_{x=1/2}=
\gamma\ln(2)-\dfrac{1}{2}\sum_{k=1}^\infty(-1)^k\dfrac{B_{2k}}{(2k)!k(2k-1)}\pi^{2k}
-\sum_{k=1}^\infty\dfrac{1}{k(2k+1)!}\psi(2k,1)\dfrac{1}{2^{2k+1}}
\]
Substraction gives
\[
\left[\int\dfrac{\ln\Gamma(x)}{x^2}dx-\int\dfrac{\ln\Gamma(x)}{(1-x)^2}dx\right]_{x=1/2}=
2-2\ln(2)-\sum_{k=1}^\infty(-1)^k\dfrac{B_{2k}}{(2k)!k(2k-1)}\pi^{2k}
\]
Substitution in \eqref{3.4aa} gives at last
\begin{equation}
S=2\ln\left(\dfrac{\pi}{2}\right)-\sum_{k=1}^\infty(-1)^k\dfrac{B_{2k}}{(2k)!k(2k-1)}\pi^{2k}
\label{3.8}
\end{equation}
For the original integral we get
\[
\int_0^1\dfrac{1}{\pi^2+\ln(x)^2}\dfrac{1}{1+x}dx=
\dfrac{1}{4}-\dfrac{1}{\pi^2}-\dfrac{1}{\pi^2}\ln\left(\dfrac{\pi}{2}\right)+
\dfrac{1}{2}\sum_{k=1}^\infty(-1)^k\dfrac{B_{2k}}{(2k)!k(2k-1)}\pi^{2k-2}
\]
This completes the proof.

\subsection{Proof of integral 8}
In this section we treat again integral $6$
\[
\int_0^1\dfrac{1}{\pi^2+\ln(x)^2}\dfrac{1}{1+x}dx
\]
to get another result. 

It is known that it is possible to convert the Bernoulli numbers to zeta numbers. We use  formula \eqref{s2.7}. Substitution in \eqref{3.8} gives after some simplification
\[
S=2\ln\left(\dfrac{\pi}{2}\right)+2\sum_{k=1}^\infty\dfrac{1}{k(2k-1)}\dfrac{1}{2^{2k}}\zeta(2k)
\]
The first factor in the summation can be split up
\[
\dfrac{1}{k(2k-1)}=\dfrac{2}{2k-1}-\dfrac{1}{k}
\]
Application gives
\[
S=2\ln\left(\dfrac{\pi}{2}\right)+4\sum_{k=1}^\infty\dfrac{1}{2k-1}\dfrac{1}{2^{2k}}\zeta(2k)-
2\sum_{k=1}^\infty\dfrac{1}{k}\dfrac{1}{2^{2k}}\zeta(2k)
\]
The last summation is known.
\[
2\sum_{k=1}^\infty\dfrac{1}{k}\dfrac{1}{2^{2k}}\zeta(2k)=2\ln\left(\dfrac{\pi}{2}\right)
\]
Substitution gives
\[
\sum_{k=1}^\infty\dfrac{\ln(k+x)}{k(k+1)}=\sum_{k=1}^\infty\dfrac{1}{2k-1}\dfrac{1}{4^{k-1}}\zeta(2k)
\]
With the integral 6 we get at last
\[
\int_0^1\dfrac{1}{\pi^2+\ln(x)^2}\dfrac{1}{1+x}dx=\dfrac{1}{4}-\dfrac{1}{\pi^2}-\dfrac{2}{\pi^2}
\sum_{k=1}^\infty\dfrac{1}{(2k-1)}\dfrac{1}{4^k}\zeta(2k)
\]
This completes the proof.

\

\subsection{Proof of integral 9}
In this section we treat integral $9$
\[
\int_0^1\dfrac{1}{4\pi^2+\ln(x)^2}\dfrac{1}{1+x}dx
\]
We use a method originating from Schl\"omilch \cite{4}. V. Moll uses this method in \cite[section 6.8]{7}. First we prove again integral 6 as a demonstration of the method. The reason doing this, is that this method is not so often described in the literature. Next we prove integral 9.

To start with the method we transform the original integral
\[
I=\int_0^1\dfrac{1}{a^2+\ln(x)^2}\dfrac{1}{1+x}dx=
\int_0^\infty\dfrac{1}{a^2+t^2}\dfrac{e^{-t/2}}{e^{t/2}+e^{-t/2}}dt
\]
The method of Schl\"omilch starts with defining a new integral using a Feynmann trick
\begin{equation}
y(x)=\int_0^\infty\dfrac{1}{a^2+t^2}\dfrac{e^{-x\, t/2}}{e^{t/2}+e^{-t/2}}dt
\label{3.4.0}
\end{equation}
and get $I=y(1)$. We note that the Feynman trick was already known to Schl\"omilch. Differentiating twice to the variable $x$ gives
\[
\dfrac{d^2 y}{dx^2}=\dfrac{1}{4}\int_0^\infty\dfrac{t^2}{a^2+t^2}\dfrac{e^{-x\, t/2}}{e^{t/2}+e^{-t/2}}dtS
\]
Multiplying the original integral with a factor $a^2$ and adding this product by the second derivative gives
\begin{equation}
4\dfrac{d^2 y}{dx^2}+a^2\, y=\int_0^\infty\dfrac{e^{-x\, t/2}}{e^{t/2}+e^{-t/2}}dt
\label{3.4.1}
\end{equation}
The integral is well-known
\[
\dfrac{1}{4}\int_0^\infty\dfrac{e^{-x\, t/2}}{e^{t/2}+e^{-t/2}}dt=
\dfrac{1}{8}\psi\left(\dfrac{3+x}{4}\right)-\dfrac{1}{8}\psi\left(\dfrac{1+x}{4}\right)=f(x)
\]
Suppose a solution of the differential equation \eqref{3.4.1} is
\[
y(x)=z_1\cos(q\pi\, x)+z_2\sin(q\pi\, x)
\]
Differentiation gives
\[
\dfrac{dy}{dx}=-z_1q\pi\sin(q\pi\, x)+z_2q\pi\cos(q\pi\, x)+
\dfrac{dz_1}{dx}cos(q\pi\, x)+\dfrac{dz_2}{dx}\sin(q\pi\, x)
\]
Suppose further
\begin{equation}
\dfrac{dz_1}{dx}\cos(q\pi\, x)+\dfrac{dz_2}{dx}\sin(q\pi\, x)=0
\label{3.4.2a}
\end{equation}
There remains
\[
\dfrac{dy}{dx}=-z_1\, q\pi\sin(q\pi\, x)+z_2\, q\pi\cos(q\pi\, x)
\]
Differentiation again gives
\[
\dfrac{d^2y}{dx^2}=-z_1(q\pi)^2\cos(q\pi\, x)-z_2(q\pi)^2\sin(q\pi\, x)
-q\pi\dfrac{z_1}{dx}\sin(q\pi\, x)+q\pi\dfrac{dz_2}{dx}\cos(q\pi\, x)
\]
Adding $\dfrac{a^2}{4}\, y$ gives
\begin{align*}
\dfrac{d^2y}{dx^2}+\dfrac{a^2}{4}\, y
&=z_1\left(\dfrac{a^2}{4}-\pi^2q^2\right)\cos(q\pi\, x)+z_2\left(\dfrac{a^2}{4}-\pi^2q^2\right)\sin(q\pi\,x) \\
&-q\pi\dfrac{z_1}{dx}\sin(q\pi\, x)+q\pi\dfrac{dz_2}{dx}\cos(q\pi\, x)
\end{align*}
Setting $q=\dfrac{a}{2\pi}$ results in
\[
\dfrac{d^2y}{dx^2}+\dfrac{a^2}{4}\, y=-q\pi\dfrac{z_1}{dx}\sin(q\pi\, x)+q\pi\dfrac{dz_2}{dx}\cos(q\pi\, x)
=f(x)
\]

Rewriting gives
\[
-\dfrac{dz_1}{dx}\sin(q\pi\, x)+\dfrac{dz_2}{dx}\cos(q\pi\, x)=\dfrac{1}{q\pi}f(x)=g(x)
\]
with
\[
g(x)=\dfrac{1}{4a}\psi\left(\dfrac{3+x}{4}\right)-\dfrac{1}{4a}\psi\left(\dfrac{1+x}{4}\right)
\]
Equation \eqref{3.4.2a} gives
\[
\dfrac{dz_1}{dx}cos(q\pi\, x)+\dfrac{dz_2}{dx}\sin(q\pi\, x)=0
\]
These two equations gives with  $q=\dfrac{a}{2\pi}$
\[
\dfrac{dz_1}{dx}=-g(x)\sin\left(\dfrac{a}{2}x\right) \qquad\qquad
\dfrac{dz_2}{dx}=g(x)\, \cos\left(\dfrac{a}{2}x\right)
\]
As solutions we try
\begin{equation}
z_1=A-\dfrac{1}{4a}\int\left(\psi\left(\dfrac{3+x}{4}\right)-\psi\left(\dfrac{1+x}{4}\right)\right)\sin\left(\dfrac{a}{2}x\right)dx
\label{3.4.3a}
\end{equation}
\begin{equation}
z_2=B+\dfrac{1}{4a}\int\left(\psi\left(\dfrac{3+x}{4}\right)-\psi\left(\dfrac{1+x}{4}\right)\right)\cos\left(\dfrac{a}{2}x\right)dx
\label{3.4.3b}
\end{equation}
With
\[
y(x)=z_1\cos\left(\dfrac{a}{2}x\right)+z_2\sin\left(\dfrac{a}{2}x\right)
\]
we get
\begin{align}
y(x)&=A\cos\left(\dfrac{a}{2}x\right)-\dfrac{1}{4a}\cos\left(\dfrac{a}{2}x\right)\int\left(\psi\left(\dfrac{3+x}{4}\right)-\psi\left(\dfrac{1+x}{4}\right)\right)\sin\left(\dfrac{a}{2}x\right)dx+ \nonumber \\
&+B\sin\left(\dfrac{a}{2}x\right)+\dfrac{1}{4a}\sin\left(\dfrac{a}{2}x\right)\int\left(\psi\left(\dfrac{3+x}{4}\right)-\psi\left(\dfrac{1+x}{4}\right)\right)\cos\left(\dfrac{a}{2}x\right)dx
\label{3.11aa}
\end{align}
Looking at this equation we see that there is a free parameter $a$. Choosing $a=\pi$ we get the original integral. We also do the derivation with $a=2\pi$.

\subsubsection{The case $a=\pi$}

For $a=\pi$ we can derive the constants $A$ and $B$ with the boundary conditions for $x=1$ and $x=-1$. Substitution in \eqref{3.11aa} gives
\begin{align}
y(x)&=A\cos\left(\dfrac{\pi}{2}\, x\right)-\dfrac{1}{4\pi}\cos\left(\dfrac{\pi}{2}\, x\right)
\int\left(\psi\left(\dfrac{3+x}{4}\right)-\psi\left(\dfrac{1+x}{4}\right)\right)\sin\left(\dfrac{\pi}{2}\, x\right)dx+ \nonumber \\
&+B\sin\left(\dfrac{\pi}{2}\, x\right)+\dfrac{1}{4\pi}\sin\left(\dfrac{\pi}{2}\, x\right)
\int\left(\psi\left(\dfrac{3+x}{4}\right)-\psi\left(\dfrac{1+x}{4}\right)\right)\cos\left(\dfrac{\pi}{2}\, x\right)dx
\label{3.4.c}
\end{align}

\

For $x=1$ equation \eqref{3.4.c}  gives
\[
y(1)=B+\dfrac{1}{4\pi}
\left[\int\left(\psi\left(\dfrac{3+x}{4}\right)-\psi\left(\dfrac{1+x}{4}\right)\right)\cos\left(\dfrac{\pi}{2}\, x\right)dx\right]_{x=1}
\]
For $x=-1$ equation \eqref{3.4.c}  gives
\[
y(-1)=-B-\dfrac{1}{4\pi}
\left[\int\left(\psi\left(\dfrac{3+x}{4}\right)-\psi\left(\dfrac{1+x}{4}\right)\right)\cos\left(\dfrac{\pi}{2}\, x\right)dx\right]_{x=1}
\]
Addition of these equations gives
\begin{equation}
y(1)=-y(-1)+\dfrac{1}{4\pi}
\int_{-1}^1\left(\psi\left(\dfrac{3+x}{4}\right)-\psi\left(\dfrac{1+x}{4}\right)\right)\cos\left(\dfrac{\pi}{2}\, x\right)dx
\label{3.10a}
\end{equation}
Equation \eqref{3.4.0} gives
\[
y(-1)=\int_0^\infty\dfrac{1}{\pi^2+t^2}\dfrac{\exp(t/2)}{\exp(-t/2)+\exp(t)}dt
\]
For the second fraction we can write
\[
\dfrac{\exp(t/2)}{\exp(-t/2)+\exp(t)}=
\dfrac{1}{2}+\dfrac{1}{2}\dfrac{\exp(t/2)-\exp(-t/2)}{\exp(t/2)+\exp(-t/2)}
\]
Substitution gives
\[
y(-1)=\dfrac{1}{2}\int_0^\infty\dfrac{1}{\pi^2+t^2}dt+
\dfrac{1}{2}\int_0^\infty\dfrac{1}{\pi^2+t^2}\dfrac{\exp(t/2)-\exp(-t/2)}{\exp(t/2)+\exp(-t/2)}dt
\]
The first integral is well-known
\[
\dfrac{1}{2}\int_0^\infty\dfrac{1}{\pi^2+t^2}dt=\dfrac{1}{2}
\]
We get
\[
y(-1)=\dfrac{1}{4}+\dfrac{1}{2}\int_0^\infty\dfrac{1}{\pi^2+t^2}\dfrac{\exp(t/2)-\exp(-t/2)}{\exp(t/2)+\exp(-t/2)}dt
\]
For the fraction in the integral we use \eqref{3.4.00} and get
\[
y(-1)=\dfrac{1}{4}+\dfrac{1}{2}\int_0^\infty\dfrac{1}{\pi^2+t^2}8\dfrac{t}{2}
\sum_{k=0}^\infty\dfrac{1}{(2k+1)^2\pi^2+t^2}dt
\]
Interchanging the integral and the summation (which is easy shown that it is allowed) gives
\[
y(-1)=\dfrac{1}{4}+2\sum_{k=0}^\infty\int_0^\infty\dfrac{t}{\pi^2+t^2}
\dfrac{1}{(2k+1)^2\pi^2+t^2}dt
\]
The integral is known and we get
\[
y(-1)+\dfrac{1}{4}+2\sum_{k=0}^\infty\dfrac{2\ln(2k+1)}{8k(k+1)\pi^2}
\]
For $k=0$ we have to take the limit. The result is
\begin{equation}
y(-1)=\dfrac{1}{4}+\dfrac{1}{\pi^2}+\dfrac{1}{2\pi^2}\sum_{k=1}^\infty\dfrac{\ln(2k+1)}{k(k+1)}
\label{3.11}
\end{equation}
Rest the integral of \eqref{3.10a}
\begin{align*}
I&=\dfrac{1}{4\pi }\int_{-1}^{1}\left( \psi \left( \dfrac{3+x}{4}\right) 
-\psi \left( \dfrac{1+x}{4}\right) \right) \cos \left( \dfrac{\pi }{2}x\right) dx= \\
&=\dfrac{1}{4\pi }\int_{0}^{1}\left( \psi \left( \dfrac{3+x}{4}\right)
-\psi \left( \dfrac{1+x}{4}\right) \right) \cos \left( \dfrac{\pi }{2}x\right) dx+ \\
&\qquad\qquad\qquad\qquad\qquad\qquad\qquad\qquad +\dfrac{1}{4\pi }\int_{-1}^{0}\left( \psi \left( \dfrac{3+x}{4}\right)-\psi \left( \dfrac{1+x}{4}\right) \right) \cos \left( \dfrac{\pi }{2}x\right) dx \\
&=\dfrac{1}{4\pi }\int_{0}^{1}\left( \psi \left( \dfrac{3+x}{4}\right)
-\psi \left( \dfrac{1+x}{4}\right) \right) \cos \left( \dfrac{\pi }{2}x\right) dx+ \\
&\qquad\qquad\qquad\qquad\qquad\qquad\qquad\qquad 
+\dfrac{1}{4\pi }\int_{0}^{1}\left( \psi \left( \dfrac{3-x}{4}\right)
-\psi \left( \dfrac{1-x}{4}\right) \right) \cos \left( \dfrac{\pi }{2}x\right) dx 
\end{align*}
\begin{align*}
&=\dfrac{1}{4\pi }\int_{0}^{1}\left( \psi \left( \dfrac{3+x}{4}\right)
-\psi \left( \dfrac{1+x}{4}\right) +\psi \left( \dfrac{3-x}{4}\right)
-\psi \left( \dfrac{1-x}{4}\right) \right) \cos \left( \dfrac{\pi }{2}x\right) dx \\
&=\dfrac{1}{4\pi }\int_{0}^{1}\left( \psi \left( 1-\dfrac{1-x}{4}\right)
-\psi \left( \dfrac{1-x}{4}\right) -\psi \left( 1-\dfrac{3-x}{4}\right)
+\psi \left( \dfrac{3-x}{4}\right) \right) \cos \left( \dfrac{\pi }{2}x\right) dx 
\end{align*}
Application of $\psi(1-z)-\psi(z)=\pi\cot(\pi\, z)$ \cite[8.365.8]{2} gives
\begin{align}
I&=\dfrac{1}{4}\int_{0}^{1}\left( \cot \left( \dfrac{1-x}{4}\pi\right)
-\cot \left( \dfrac{3-x}{4}\pi \right) \right) \cos \left( \dfrac{\pi }{2}x\right) dx \nonumber \\
&=\dfrac{1}{4}2\int_{0}^{1}\dfrac{1}{\cos \left( \dfrac{\pi }{2}x\right) }\cos \left( \dfrac{\pi }{2}x\right) dx=\dfrac{1}{2}
\label{3.12}
\end{align}
Substitution of \eqref{3.11} and \eqref{3.12} in \eqref{3.10a} gives at last
\begin{equation}
\int_0^1\dfrac{1}{\pi^2+\ln(x)^2}\dfrac{1}{1+x}dx=
\dfrac{1}{4}-\dfrac{1}{\pi^2}-\dfrac{1}{2\pi^2}\sum_{k=1}^\infty\dfrac{\ln(2k+1)}{k(k+1)}
\label{3.16aa}
\end{equation}
which we already found in a previous section. 

\subsubsection{The case $a=2\pi$}
Setting $a=2\pi$ in equation \eqref{3.11aa} gives
\begin{align*}
y(x)&=A\cos(\pi x)-\dfrac{1}{8\pi}\cos(\pi x)
\int\left(\psi\left(\dfrac{3+x}{4}\right)-\psi\left(\dfrac{1+x}{4}\right)\right)\sin(\pi x)dx+ \nonumber \\ 
&+B\sin(\pi x)+\dfrac{1}{8\pi}\sin(\pi x)
\int\left(\psi\left(\dfrac{3+x}{4}\right)-\psi\left(\dfrac{1+x}{4}\right)\right)\cos(\pi x)dx
\end{align*}
We use the boundary conditions for $x=0$ and $x=1$.
\[
y(0)=A-\dfrac{1}{8\pi}
\left[\int\left(\psi\left(\dfrac{3+x}{4}\right)-\psi\left(\dfrac{1+x}{4}\right)\right)\sin(\pi x)dx\right]_{x=0}
\]
For $y(0)$ we get from equation \eqref{3.4.0} with \cite[Table 132(2)]{1}
\[
y(0)=\int_0^\infty\dfrac{1}{4\pi^2+t^2}\dfrac{1}{e^{t/2}+e^{-t/2}}dt=
\int_0^1\dfrac{1}{4\pi^2+\ln(x)^2}\dfrac{1}{(x+1)\sqrt{x}}dx=\dfrac{4-\pi}{8\pi}
\]
Then we get for $A$
\[
A=\dfrac{4-\pi}{8\pi}+\dfrac{1}{8\pi}
\left[\int\left(\psi\left(\dfrac{3+x}{4}\right)-\psi\left(\dfrac{1+x}{4}\right)\right)\sin(\pi x)dx\right]_{x=0}
\]
For $y(1)$ we get
\[
y(1)=-A+\dfrac{1}{8\pi}
\left[\int\left(\psi\left(\dfrac{3+x}{4}\right)-\psi\left(\dfrac{1+x}{4}\right)\right)\sin(\pi x)dx\right]_{x=1}
\]
Substitution of $A$ gives
\begin{align*}
y(1)&=-\dfrac{4-\pi}{8\pi}-\dfrac{1}{8\pi}
\left[\int\left(\psi\left(\dfrac{3+x}{4}\right)-\psi\left(\dfrac{1+x}{4}\right)\right)\sin(\pi x)dx\right]_{x=0}+ \\
&\qquad\quad\quad\ \ +\dfrac{1}{8\pi}
\left[\int\left(\psi\left(\dfrac{3+x}{4}\right)-\psi\left(\dfrac{1+x}{4}\right)\right)\sin(\pi x)dx\right]_{x=1}
\end{align*}
Rewriting gives
\[
y(1)=-\dfrac{4-\pi}{8\pi}+\dfrac{1}{8\pi}
\int_0^1\left(\psi\left(\dfrac{3+x}{4}\right)-\psi\left(\dfrac{1+x}{4}\right)\right)\sin(\pi x)dx
\]
The integral can be computed with partial integration. The result is
\[
y(1)=\dfrac{1}{8}-\dfrac{1}{2\pi}+2\int_{3/4}^1\ln\Gamma(y)\cos(4\pi y)dy-
2\int_{1/4}^{2/4}\ln\Gamma(y)\cos(4\pi y)dy
\]
For the $\ln\Gamma(y)$ function we use the well-known series representation due to Kummer \cite[1.91.(14)]{6}
\[
\ln\Gamma(y)=\left(\dfrac{1}{2}-y\right)\big(\gamma+\ln(2)\big)+(1-y)\ln(\pi)-\dfrac{1}{2}\ln\sin(\pi y)+
\sum_{k=1}^\infty\dfrac{\ln(k)}{k}\sin(2\pi k y)
\]
Substitution in the integrals gives
\[
y(1)=\dfrac{1}{8}+2\int_{3/4}^1\sum_{k=1}^\infty\dfrac{\ln(k)}{k}\sin(2\pi k y)\cos(4\pi y)dy-
2\int_{1/4}^{2/4}\sum_{k=1}^\infty\dfrac{\ln(k)}{k}\sin(2\pi k y)\cos(4\pi y)dy
\]
After interchanging the integrals and the summations (which is allowed) we get
\[
y(1)=\dfrac{1}{8}+2\sum_{k=1}^\infty\dfrac{\ln(k)}{k}\left(\int_{3/4}^1\sin(2\pi k y)\cos(4\pi y)dy-
\int_{1/4}^{2/4}\sin(2\pi k y)\cos(4\pi y)dy\right)
\]
The integrals are elementary. 
\[
\int_{3/4}^1\sin(2\pi k y)\cos(4\pi y)dy-\int_{1/4}^{2/4}\sin(2\pi k y)\cos(4\pi y)dy=\cos(k \pi)-1=(-1)^k-1
\]
The result for the original integral is at last
\[
\int_0^1\dfrac{1}{4\pi^2+\ln(x)^2}\dfrac{1}{1+x}dx=
\dfrac{1}{8}-\dfrac{2}{\pi^2}\sum_{k=1}^\infty\dfrac{\ln(2k+1)}{(2k-1)(2k+3)}
\]
This completes the proof.

\subsection{Proof of integral 10}
In this section we treat integral $10$
\[
I=\int_0^1\dfrac{1}{\left(a^2+\ln(x)^2\right)^2}\dfrac{1}{1+x}dx
\]
We start with integral 3
\begin{equation}
\int_0^1\dfrac{1}{a^2+\ln(x)^2}\dfrac{1}{1+x}dx=
\dfrac{\pi}{4a}-2\sum_{k=0}^\infty\dfrac{\ln(a)-\ln(\pi)-\ln(2k+1)}{a^2-(2k+1)^2\pi^2}
\label{3.20}
\end{equation}
Differentiation of the left hand side of \eqref{3.20} to the parameter $a$ gives
\begin{equation}
\dfrac{d}{da}\int_0^1\dfrac{1}{a^2+\ln(x)^2}\dfrac{1}{1+x}dx=
-2a\int_0^1\dfrac{1}{\left(a^2+\ln(x)^2\right)^2}\dfrac{1}{1+x}dx
\label{3.21}
\end{equation}
Differentiation of the right hand side of \eqref{3.20} with respect to the parameter $a$ gives
\begin{multline}
\dfrac{d}{da}\left(\dfrac{\pi}{4a}-2\sum_{k=0}^\infty\dfrac{\ln(a)-\ln(\pi)-\ln(2k+1)}{a^2-(2k+1)^2\pi^2}\right)=  \\
=-\dfrac{\pi}{4a^2}-\dfrac{2}{a}\sum_{k=0}^\infty\dfrac{1}{a^2-(\pi-2k\pi)^2}+
4a\sum_{k=0}^\infty\dfrac{\ln(\pi+2k\pi)-\ln(a)}{\big((\pi+2k\pi)^2-a^2\big)^2}
\label{3.23}
\end{multline}
The first summation of the right hand side is known. We get
\begin{multline}
\dfrac{d}{da}\left(\dfrac{\pi}{4a}-2\sum_{k=0}^\infty\dfrac{\ln(a)-\ln(\pi)-\ln(2k+1)}{a^2-(2k+1)^2\pi^2}\right)=  \\
=-\dfrac{\pi}{4a^2}+\dfrac{1}{2a^2}\tan\left(\dfrac{a}{2}\right)+
4a\sum_{k=0}^\infty\dfrac{\ln(\pi+2k\pi)-\ln(a)}{\big((\pi-2k\pi)^2-a^2\big)^2}
\label{3.22}
\end{multline}
Combination of \eqref{3.21} and \eqref{3.22} gives
\[
-2a\int_0^1\dfrac{1}{\left(a^2+\ln(x)^2\right)^2}\dfrac{1}{1+x}dx=
-\dfrac{\pi}{4a^2}+\dfrac{1}{2a^2}\tan\left(\dfrac{a}{2}\right)-
4a\sum_{k=0}^\infty\dfrac{\ln(\pi+2k\pi)-\ln(a)}{\big((\pi-2k\pi)^2-a^2\big)^2}
\]
Simplification gives the desired result
\begin{equation}
\int_0^1\dfrac{1}{\left(a^2+\ln(x)^2\right)^2}\dfrac{1}{1+x}dx=
\dfrac{\pi}{8a^3}-\dfrac{1}{4a^3}\tan\left(\dfrac{a}{2}\right)+
2\sum_{k=0}^\infty\dfrac{\ln(\pi+2k\pi)-\ln(a)}{\big((\pi-2k\pi)^2-a^2\big)^2}
\label{3.24}
\end{equation}
This completes the proof.

\subsection{Proof of integral 11}
In this section we treat integral $11$
\[
\int_0^1\dfrac{1}{\left(\pi^2+\ln(x)^2\right)^2}\dfrac{1}{1+x}dx
\]
When taking the limit for $a \rightarrow \pi$ in equation \eqref{3.24} there arises a problem for $k=0$. So we had to split up equation \eqref{3.23} and get
\begin{align*}
&\left[\dfrac{d}{da}\left(\dfrac{\pi}{4a}-2\sum_{k=0}^\infty\dfrac{\ln(a)-\ln(\pi)-\ln(2k+1)}{a^2-(2k+1)^2\pi^2}\right)\right]_{a=\pi}=  \\
&=\left[-\dfrac{\pi}{4a^2}-\dfrac{2}{a}\dfrac{1}{a^2-\pi^2}+
4a\dfrac{\ln(a)-\ln(\pi)}{(a^2-\pi^2)^2}-
\dfrac{2}{a}\sum_{k=1}^\infty\dfrac{1}{a^2-(\pi+2k\pi)^2}+
4a\sum_{k=1}^\infty\dfrac{\ln(a)-\ln(\pi+2k\pi)}{\big((\pi+2k\pi)^2-a^2\big)^2}\right]_{a=\pi}=
\end{align*}
\begin{align*}
&=-\dfrac{1}{4\pi}-\lim_{a \rightarrow \pi}\left(\dfrac{2}{a}\dfrac{1}{a^2-\pi^2}-
4a\dfrac{\ln(a)-\ln(\pi)}{(a^2-\pi^2)^2}\right)+\dfrac{1}{2\pi^3}-
\dfrac{1}{4\pi^3}\sum_{k=1}^\infty\dfrac{\ln(2k+1)}{\big(k(1+k)\big)^2}
\qquad\qquad\qquad\qquad \\
&=-\dfrac{1}{4\pi}+\dfrac{1}{\pi^3}+\dfrac{1}{2\pi^3}-
\dfrac{1}{4\pi^3}\sum_{k=1}^\infty\dfrac{\ln(2k+1)}{\big(k(1+k)\big)^2}
=-\dfrac{1}{4\pi}+\dfrac{3}{2\pi^3}-
\dfrac{1}{4\pi^3}\sum_{k=1}^\infty\dfrac{\ln(2k+1)}{\big(k(1+k)\big)^2} 
\end{align*}
Combination of \eqref{3.21} with $a=\pi$ and the last equation gives
\[
-2\pi\int_0^1\dfrac{1}{\left(\pi^2+\ln(x)^2\right)^2}\dfrac{1}{1+x}dx=
-\dfrac{1}{4\pi}+\dfrac{3}{2\pi^3}-
\dfrac{1}{4\pi^3}\sum_{k=1}^\infty\dfrac{\ln(1+2k)}{\big(k(1+k)\big)^2}
\]
Simplification gives
\[
\int_0^1\dfrac{1}{\left(\pi^2+\ln(x)^2\right)^2}\dfrac{1}{1+x}dx=
\dfrac{1}{8\pi^2}-\dfrac{3}{4\pi^4}+\dfrac{1}{8\pi^4}\sum_{k=1}^\infty\dfrac{\ln(1+2k)}{k^2(1+k)^2}
\]
This completes the proof.

\subsection{Proof of integral 12}
In this section we treat integral $12$
\[
I=\int_0^1\dfrac{1}{a^2+\ln(x)^2}\dfrac{1}{(1+x)^2}dx
\]
and start the proof with partial integration. 
\begin{align*}
&u=\dfrac{1}{a^2+\ln(x)^2}	\qquad\qquad u'=-\dfrac{2\ln)x)}{x(a^2+\ln(x)^2)^2} \\
&v'=\dfrac{1}{(1+x)^2} 	\qquad\qquad\quad v=-\dfrac{1}{1+x}
\end{align*}
Application gives
\[
I=\left[-\dfrac{1}{a^2+\ln(x)^2}\dfrac{1}{1+x}\right]_0^1-
\int_0^1\dfrac{2\ln)x)}{(a^2+\ln(x)^2)^2}\dfrac{1}{x(1+x)}dx
\]
Some manipulation gives
\[
I=-\dfrac{1}{2a^2}-2\int_0^1\dfrac{\ln)x)}{(a^2+\ln(x)^2)^2}\dfrac{1}{x}dx+
2\int_0^1\dfrac{\ln)x)}{(a^2+\ln(x)^2)^2}\dfrac{1}{1+x}dx
\]
The first integral is elementary.
\[
2\int_0^1\dfrac{\ln)x)}{(a^2+\ln(x)^2)^2}\dfrac{1}{x}dx=-\dfrac{1}{a^2}
\]
Substitution gives
\[
I=\dfrac{1}{2a^2}+2\int_0^1\dfrac{\ln)x)}{(a^2+\ln(x)^2)^2}\dfrac{1}{1+x}dx
\]
Setting $b=a^2$ and integrating to the parameter $b $ gives
\begin{align*}
\int I(b)db
&=\int\dfrac{1}{2b}db+2\int_0^1\left(\int\dfrac{1}{(b+\ln(x)^2)^2}db\right)\dfrac{\ln(x)}{1+x}dx \\
&=\dfrac{1}{2}\ln(b)-2\int_0^1 \dfrac{\ln(x)}{b+\ln(x)^2}\dfrac{1}{1+x}dx
\end{align*}
The last integral is the integral in nr. 14 in table 129 of Bierens de Haan (see Appendix)
\[
\int I(b)db=\dfrac{1}{2}\ln(b)-\ln\left(\dfrac{2\sqrt{b}}{\pi}\right)-
\psi\left(\dfrac{\sqrt{b}}{2\pi}\right)+2\psi\left(\dfrac{\sqrt{b}}{\pi}\right)
\]
Differentiation to the parameter $b$ gives
\[
I(b)=\dfrac{1}{2b}-\dfrac{1}{2b}-\dfrac{1}{4\sqrt{b}\pi}\psi^{(1)}\left(\dfrac{\sqrt{b}}{2\pi}\right)
+\dfrac{1}{\sqrt{b}\pi}\psi^{(1)}\left(\dfrac{\sqrt{b}}{\pi}\right)
\]
Substitution of $b=a^2$ gives at last
\[
\int_0^1\dfrac{1}{a^2+\ln(x)^2}\dfrac{1}{(1+x)^2}dx=-\dfrac{1}{4a\pi}\psi^{(1)}\left(\dfrac{a}{2\pi}\right)
+\dfrac{1}{a\pi}\psi^{(1)}\left(\dfrac{a}{\pi}\right)
\]
This completes the proof.

\subsection{Proof of integral 13}
In this section we prove the integrals 13 
\[
\int_0^1\dfrac{1}{\pi^2+\ln(x)^2}\dfrac{1}{(1+x)^2}dx=
\int_1^\infty\dfrac{1}{\pi^2+\ln(x)^2}\dfrac{1}{(1+x)^2}dx=
\dfrac{1}{24}
\]
We use two methods. The first method uses the result of the previous integral. The second method uses a Feynmann trick.
\subsubsection{Direct proof}
Setting $a=\pi$ in the integral nr. 12 gives for the different terms
\[
-\dfrac{1}{4\pi^2}\psi^{(1)}\left(\dfrac{1}{2}\right)=-\dfrac{1}{8}
\qquad\qquad
\dfrac{1}{\pi^2}\psi^{(1)}(1,1)=\dfrac{1}{6}
\]
This gives the first result
\[
\int_0^1\dfrac{1}{\pi^2+\ln(x)^2}\dfrac{1}{(1+x)^2}dx=\dfrac{1}{24}
\]
Transforming the first integral with $x \rightarrow \dfrac{1}{x}$ gives
\[
\int_0^1\dfrac{1}{\pi^2+\ln(x)^2}\dfrac{1}{(1+x)^2}dx=
\int_1^\infty\dfrac{1}{\pi^2+\ln(x)^2}\dfrac{1}{(1+x)^2}dx=\dfrac{1}{24}
\]
This completes the proof.

Adding the two integrals gives
\[
\int_0^\infty\dfrac{1}{\pi^2+\ln(x)^2}\dfrac{1}{(1+x)^2}dx=\dfrac{1}{12}
\]

\subsubsection{Using a Feynmann trick}
We start with the same trick as in section 4.2.2.
\[
\int_0^1\dfrac{1}{\pi^2+\ln(x)^2}\dfrac{1}{(1+x)^2}dx \rightarrow 
\int_0^1\dfrac{1}{\pi^2+\ln(x)^2}\dfrac{1}{(1+b\, x)^2}dx
\]
After integration and differentiation twice to the parameter $b$ we get after interchanging the integral and the summation
\[
\int_0^1\dfrac{1}{\pi^2+\ln(x)^2}\dfrac{1}{(1+b\, x)^2}dx
=-\dfrac{d^2}{db^2}\int_0^1\dfrac{1}{\pi^2+\ln(x)^2}\dfrac{\ln(1+b\, x)}{x^2}dx \\
\]
The last fraction can be written as a power series
\begin{align*}
\int_0^1\dfrac{1}{\pi^2+\ln(x)^2}\dfrac{1}{(1+b\, x)^2}dx
&=-\dfrac{d^2}{db^2}\int_0^1\dfrac{1}{\pi^2+\ln(x)^2}
\sum_{k=1}^\infty(-1)^{k-1}\dfrac{b^k}{k}x^{k-2}dx \nonumber \\
&=\dfrac{d^2}{db^2}\sum_{k=1}^\infty(-1)^k\dfrac{b^k}{k}
\int_0^1\dfrac{1}{\pi^2+\ln(x)^2}x^{k-2}dx \nonumber \\
&=\sum_{k=2}^\infty(-1)^k(k-1)b^{k-2}\int_0^1\dfrac{1}{\pi^2+\ln(x)^2}x^{k-2}dx
\end{align*}
The integral is known.
\[
\int_0^1\dfrac{1}{\pi^2+\ln(x)^2}x^{k-2}dx=\dfrac{(-1)^k}{\pi}
\left(\dfrac{i}{2}\big(\text{Ei}(-i(k-1)\pi)-\text{Ei}(i(k-1)\pi)\big)-\pi\right)
\]
Application of \eqref{s2.6} gives
\[
\int_0^1\dfrac{1}{\pi^2+\ln(x)^2}x^{k-2}dx=
\dfrac{(-1)^k}{\pi}\left(\text{Si}((k-1)\pi)-\dfrac{\pi}{2}\right)
\]
Setting $b=1$ gives
\[
\int_0^1\dfrac{1}{\pi^2+\ln(x)^2}\dfrac{1}{(1+x)^2}dx=
\dfrac{1}{\pi}\sum_{k=2}^\infty(k-1)\left(\text{Si}((k-1)\pi)-\dfrac{\pi}{2}\right)
\]
Rewriting the summation from $k=1$ to $\infty$ gives
\[
\int_0^1\dfrac{1}{\pi^2+\ln(x)^2}\dfrac{1}{(1+x)^2}dx=
\dfrac{1}{\pi}\sum_{k=1}^\infty k\left(\text{Si}(k\pi)-\dfrac{\pi}{2}\right)=
\dfrac{1}{\pi}\sum_{k=1}^\infty k\, \text{si}(k\pi)
\]
This completes the proof.

\subsection{Proof of integral 14}
In this section we treat integral $14$
\[
\int_0^1\dfrac{\ln(x)}{a^2+\ln(x)^2}\dfrac{1}{1+x}dx
\]
We use two methods. The first method uses integral nr. 14 of the list of Bierens de Haan (see Appendix). The second method uses Legendre's formula.

\subsubsection{Using nr. 14 of the list in the Appendix}

We start the proof with integral nr. 14 of the Appendix.
\[
\int_0^1\dfrac{\ln(x)}{q^2+\ln(x)^2}\dfrac{x}{1-x^2}dx=
\dfrac{1}{2}\left(\dfrac{\pi}{2q}+\ln\left(\dfrac{\pi}{q}\right)+\psi\left(\dfrac{q}{\pi}\right)\right)
\]
We split the second factor in the integral
\[
\dfrac{x}{1-x^2}=\dfrac{1}{2(1-x)}-\dfrac{1}{2(1+x)}
\]
Substitution with $q=a$ gives
\[
\int_0^1\dfrac{\ln(x)}{a^2+\ln(x)^2}\dfrac{1}{1-x}dx-
\int_0^1\dfrac{\ln(x)}{a^2+\ln(x)^2}\dfrac{1}{1+x}dx=
\dfrac{\pi}{2a}+\ln\left(\dfrac{\pi}{a}\right)+\psi\left(\dfrac{a}{\pi}\right)
\]
The first integral is known as nr. 2 of the list of Bierens de Haan.
\[
\int_0^1\dfrac{\ln(x)}{q^2+ \ln(x)^2}\dfrac{1}{1-x}dx=
\dfrac{1}{2}\left(\dfrac{\pi}{q}+\ln\left(\dfrac{2\pi}{q}\right)+\psi\left(\dfrac{q}{2\pi}\right)\right) 
\]
Substitution gives
\[
\int_0^1\dfrac{\ln(x)}{a^2+\ln(x)^2}\dfrac{1}{1+x}dx=
\dfrac{1}{2}\left(\dfrac{\pi}{a}+\ln\left(\dfrac{2\pi}{a}\right)+\psi\left(\dfrac{a}{2\pi}\right)\right)-
\dfrac{\pi}{2a}-\ln\left(\dfrac{\pi}{a}\right)-\psi\left(\dfrac{a}{\pi}\right)
\]
Using 
\[
\dfrac{1}{2}\ln\left(\dfrac{2\pi}{a}\right)-\ln\left(\dfrac{\pi}{a}\right)=\dfrac{1}{2}\ln\left(\dfrac{2a}{\pi}\right)
\]
gives
\begin{equation}
\int_0^1\dfrac{\ln(x)}{a^2+\ln(x)^2}\dfrac{1}{1+x}dx=
\dfrac{1}{2}\ln\left(\dfrac{2a}{\pi}\right)+
\dfrac{1}{2}\psi\left(\dfrac{a}{2\pi}\right)-\psi\left(\dfrac{a}{\pi}\right)
\label{4.23bb}
\end{equation}
and this completes the proof.

Note: We checked that the nrs. 2 and 14 from the list of Bierens de Haan are correct.

\subsubsection{Using Legendre's formula}
We use Legendre's formula \cite[Exerc. 40 b.]{3}
\begin{equation}
\int_0^\infty\dfrac{t}{(e^{b\, t}+1)(t^2+a^2)}dt=
\dfrac{1}{2}\psi\left(\dfrac{1}{2}+\dfrac{a\, b}{2\pi}\right)-\dfrac{1}{2}\ln\left(\dfrac{a\, b}{2\pi}\right)
\label{4.23aa}
\end{equation}
Applying the substitution $t=-\ln(x)$ with $b=1$ gives
\[
\int_0^\infty\dfrac{t}{(e^{t}+1)(t^2+a^2)}dt=
-\int_0^1\dfrac{\ln(x)}{a^2+\ln(x)^2}\dfrac{1}{1+x}dx
\]
Using \eqref{4.23aa} gives
\[
\int_0^\infty\dfrac{t}{(e^{t}+1)(t^2+a^2)}dt=
\dfrac{1}{2}\ln\left(\dfrac{a}{2\pi}\right)-\dfrac{1}{2}\psi\left(\dfrac{1}{2}+\dfrac{a}{2\pi}
\right)
\]
For the function $\psi(x)$ we can use
\[
\psi\left(\dfrac{a}{\pi}\right)=\dfrac{1}{2}\psi\left(\dfrac{1}{2}+\dfrac{a}{2\pi}\right)+\dfrac{1}{2}\psi\left(\dfrac{a}{2\pi}\right)+\ln(2)
\]
and get the same result as \eqref{4.23bb}. This completes the proof.

\subsection{Proof of integral 15}
In this section we treat integral $15$
\[
\int_0^1\dfrac{\ln(x)}{\pi^2+\ln(x)^2}\dfrac{1}{1+x}dx
\]
We  use two methods. The first method uses the result of the previous integral. The second method uses power series.
\subsubsection{Direct proof}
Setting $a=\pi$ in the integral of 14 and using
\[
\psi(1)=-\gamma \qquad \psi\left(\dfrac{1}{2}\right)=-\gamma-2\ln(2)
\]
gives
 \[
\int_0^1\dfrac{\ln(x)}{\pi^2+\ln(x)^2}\dfrac{1}{1+x}dx=\dfrac{\gamma}{2}-\dfrac{1}{2}\ln(2)
\]
and this completes the proof.

\subsubsection{Using power series}
Because $x$ lies between $0$ and $1$ the second factor in the integral can be written as a geometric series
\[
\dfrac{1}{1+x}=\sum_{k=0}^\infty(-x)^k 
\]
The integral becomes
\[
I=\int_0^1\dfrac{\ln(x)}{a^2+\ln(x)^2}\dfrac{1}{1+x}dx=\int_0^1\dfrac{\ln(x)}{a^2+\ln(x)^2}\sum_0^\infty(-x)^kdx
\]
Because both factors are absolute convergent we may interchange the summation and the integral
\[
I=\int_0^1\dfrac{\ln(x)}{a^2+\ln(x)^2}\dfrac{1}{1+x}dx=
\sum_{k=0}^\infty\int_0^1\dfrac{\ln(x)(-x)^k}{a^2+\ln(x)^2}dx
\]
The integral is known and we get
\begin{multline*}
\int_0^1\dfrac{\ln(x)(-x)^k}{a^2+\ln(x)^2}dx= \\
=\dfrac{1}{2}(-1)^k\Big(2\cos\big(a(1+k)\big)\text{Ci}\big(a(1+k)\big)-\sin\big(a(1+k)\big)(\pi-2\text{Si}\big(a(1+k)\big)\Big)
\end{multline*}

To simplify this integral we set $a=\pi$ and get
\[
\int_0^1\dfrac{\ln(x)(-x)^k}{\pi^2+\ln(x)^2}dx=-\text{Ci}\big((k+1)\pi\big)
\]
Rewriting the summation gives
\[
\int_0^1\dfrac{\ln(x)}{\pi^2+\ln(x)^2}\dfrac{1}{1+x}dx=-\sum_{k=1}^\infty\text{Ci}(k\, \pi)
\]
Combining both methods results in
\[
\int_0^1\dfrac{\ln(x)}{\pi^2+\ln(x)^2}\dfrac{1}{1+x}dx=\sum_{k=1}^\infty\text{Ci}(k\, \pi)=\dfrac{1}{2}\ln(2)-\dfrac{\gamma}{2}
\]
and this completes the proof.

\subsection{Proof of integral 16}
In this section we treat integral $16$
\[
\int_0^1\dfrac{\ln(x)}{\left(a^2+\ln(x)^2\right)^2}\dfrac{1}{1+x}dx
\]
We  use two methods. The first method uses differentiation. The second method uses partial integration.

\subsubsection{Differentiation}
Setting $a^2=b$ in the integral 14 and differentiating both sides of the integral with respect to the parameter $b$ gives
\[
\dfrac{d}{db}\int_0^1\dfrac{\ln(x)}{b+\ln(x)^2}\dfrac{1}{1+x}dx=
-\int_0^1\dfrac{\ln(x)}{(b+\ln(x)^2)^2}\dfrac{1}{1+x}dx \\
\]
Replacing $b=a^2$ gives
\[
\int_0^1\dfrac{\ln(x)}{(a^2+\ln(x)^2)^2}\dfrac{1}{1+x}dx=
-\dfrac{1}{4a^2}-\dfrac{1}{8a\pi}\psi^{(1)}\left(\dfrac{a}{2\pi}\right)+\dfrac{1}{2a\pi}\psi^{(1)}\left(\dfrac{a}{\pi}\right)
\]
This completes the proof.
\subsubsection{Partial integration}
We start with integral 12
\[
\int_0^1\dfrac{1}{a^2+\ln(x)^2}\dfrac{1}{(1+x)^2}dx=
-\dfrac{1}{4a\pi}\psi^{(1)}\left(\dfrac{a}{2\pi}\right)+\dfrac{1}{a\pi}\psi^{(1)}\left(\dfrac{a}{\pi}\right)
\]
Setting $u=\dfrac{1}{a^2+\ln(x)^2}$ and $v'=\dfrac{1}{(1+x)^2}$ gives
\begin{align*}
\int_0^1\dfrac{1}{a^2+\ln(x)^2}\dfrac{1}{(1+x)^2}dx
&=\left[-\dfrac{1}{1+x}\dfrac{1}{a^2+\ln(x)^2}\right]_0^1-\int_0^1\dfrac{2\ln(x)}{x\big(a^2+\ln(x)^2\big)^2}\dfrac{1}{1+x}dx \\
&=-\dfrac{1}{2a^2}-2\int_0^1\dfrac{\ln(x)}{\big(a^2+\ln(x)^2\big)^2}\dfrac{1}{x}dx+2\int_0^1\dfrac{\ln(x)}{\big(a^2+\ln(x)^2\big)^2}\dfrac{1}{1+x}dx \\
&=\dfrac{1}{2a^2}+2\int_0^1\dfrac{\ln(x)}{\big(a^2+\ln(x)^2\big)^2}\dfrac{1}{1+x}dx
\end{align*}
Rewriting gives
\[
\int_0^1\dfrac{\ln(x)}{\big(a^2+\ln(x)^2\big)^2}\dfrac{1}{1+x}dx=
\dfrac{1}{2}\int_0^1\dfrac{1}{a^2+\ln(x)^2}\dfrac{1}{(1+x)^2}dx-\dfrac{1}{4a^2}
\]
Using integral 12 gives
\[
\int_0^1\dfrac{\ln(x)}{(a^2+\ln(x)^2)^2}\dfrac{1}{1+x}dx=
-\dfrac{1}{4a^2}-\dfrac{1}{8a\pi}\psi^{(1)}\left(\dfrac{a}{2\pi}\right)+\dfrac{1}{2a\pi}\psi^{(1)}\left(\dfrac{a}{\pi}\right)
\]
This completes the proof.

\subsection{Proof of integral 17}
In this section we treat integral $17$
\[
\int_0^1\dfrac{\ln(x)}{\left(\pi^2+\ln(x)^2\right)^2}\dfrac{1}{1+x}dx
\]
Using $a=\pi$ and
\[
\psi^{(1)}\left(\dfrac{1}{2}\right)=\dfrac{\pi^2}{2} \qquad\qquad
\psi^{(1)}(1)=\dfrac{\pi^2}{6}
\]
and substituting these values in the integral 16 gives
\[
\int_0^1\dfrac{\ln(x)}{\left(\pi^2+\ln(x)^2\right)^2}\dfrac{1}{1+x}dx=
\dfrac{1}{48}-\dfrac{1}{4\pi^2}
\]
and this completes the proof

\subsection{Proof of integral 18}
In this section we treat integral $18$
\[
\int_0^1\dfrac{\ln(x)}{\pi^2+\ln(x)^2}\dfrac{1}{(1+x)^2}dx
\]
We  use two methods. The first method uses power series. The second method uses a Feynmann trick.

\subsubsection{Using power series}
Because $x$ lies between $0$ and $1$ the second factor in the integral can be written as a geometric series
\[
\dfrac{1}{(1+x)^2}=\sum_{k=0}^\infty(k+1)(-x)^k 
\]
The integral becomes
\[
I=\int_0^1\dfrac{\ln(x)}{a^2+\ln(x)^2}\dfrac{1}{(1+x)^2}dx=
\int_0^1\dfrac{\ln(x)}{a^2+\ln(x)^2}\sum_{k=0}^\infty(k+1)(-x)^kdx
\]
Because both factors are absolute convergent we may interchange the summation and the integral
\[
I=\int_0^1\dfrac{\ln(x)}{a^2+\ln(x)^2}\dfrac{1}{(1+x)^2}dx=
\sum_{k=0}^\infty(k+1)\int_0^1\dfrac{\ln(x)(-x)^k}{a^2+\ln(x)^2}dx
\]
The integral is known and we get
\begin{multline*}
\int_0^1\dfrac{\ln(x)(-x)^k}{a^2+\ln(x)^2}dx= \\
=\dfrac{1}{2}(-1)^k\Big(2\cos\big(a(1+k)\big)\text{Ci}\big(a(1+k)\big)-\sin\big(a(1+k)\big)(\pi-2\text{Si}\big(a(1+k)\big)\Big)
\end{multline*}

To simplify this integral we set $a=\pi$ and get
\[
\int_0^1\dfrac{\ln(x)(-x)^k}{\pi^2+\ln(x)^2}dx=-\text{Ci}\big((k+1)\pi\big)
\]
Rewriting the summation gives
\[
\int_0^1\dfrac{\ln(x)}{\pi^2+\ln(x)^2}\dfrac{1}{(1+x)^2}dx=
-\sum_{k=1}^\infty\, k\, \text{Ci}(k\, \pi)
\]

\subsubsection{Using a Feynmann trick}
We start with the same trick as in section 4.2.2.
\[
\int_0^1\dfrac{\ln(x)}{a^2+\ln(x)^2}\dfrac{1}{(1+x)^2}dx \rightarrow 
\int_0^1\dfrac{\ln(x)}{a^2+\ln(x)^2}\dfrac{1}{(1+b\, x)^2}dx
\]
After integration and differentiation twice to the parameter $b$ we get after interchanging the integral and the summation
\[
\int_0^1\dfrac{\ln(x)}{a^2+\ln(x)^2}\dfrac{1}{(1+b\, x)^2}dx
=-\dfrac{d^2}{db^2}\int_0^1\dfrac{\ln(x)}{a^2+\ln(x)^2}\dfrac{\ln(1+b\, x)}{x^2}dx \\
\]
The last fraction can be written as a power series
\begin{align}
\int_0^1\dfrac{\ln(x)}{a^2+\ln(x)^2}\dfrac{1}{(1+b\, x)^2}dx
&=-\dfrac{d^2}{db^2}\int_0^1\dfrac{\ln(x)}{a^2+\ln(x)^2}
\sum_{k=1}^\infty(-1)^{k-1}\dfrac{b^k}{k}x^{k-2}dx \nonumber \\
&=\dfrac{d^2}{db^2}\sum_{k=1}^\infty(-1)^k\dfrac{b^k}{k}
\int_0^1\dfrac{\ln(x)}{a^2+\ln(x)^2}x^{k-2}dx \nonumber \\
&=\sum_{k=1}^\infty(-1)^k(k-1)b^{k-2}\int_0^1\dfrac{\ln(x)}{a^2+\ln(x)^2}x^{k-2}dx
\label{3.25}
\end{align}
The integral is known
\begin{multline*}
\int_0^1\dfrac{\ln(x)}{a^2+\ln(x)^2}x^{k-2}dx= \\
=\dfrac{1}{4a}\big[e^{ia(k+1)}(-i+a+ak)\text{Ei}(-ia(1+k))+e^{-ia(k+1)}(i+a+ak)\text{Ei}(ia(1+k))\big]
\end{multline*}
To simplify this integral we set $a=\pi$ and get
\[
\int_0^1\dfrac{\ln(x)}{\pi^2+\ln(x)^2}x^{k-2}dx=
-(-1)^k\dfrac{1}{2}\Big[\text{Ei}(-i\pi(1-k))+\text{Ei}(i\pi(1-k))\Big]
\]
Substitution in \eqref{3.25} and setting $b=1$ gives
\[
\int_0^1\dfrac{\ln(x)}{\pi^2+\ln(x)^2}\dfrac{1}{(1+x)^2}dx
=-\sum_{k=1}^\infty(k-1)\dfrac{1}{2}\Big[\text{Ei}(-i\pi(1-k))+\text{Ei}(i\pi(1-k))\Big]
\]
Because for $k=1$ the first term is zero we rewrite the summation as
\[
\int_0^1\dfrac{\ln(x)}{\pi^2+\ln(x)^2}\dfrac{1}{(1+x)^2}dx
=-\sum_{k=1}^\infty\, k\dfrac{1}{2}\Big[\text{Ei}(-i\pi\, k)+\text{Ei}(i\pi\, k)\Big]
\]
Using \eqref{s2.6} gives
\[
\int_0^1\dfrac{\ln(x)}{\pi^2+\ln(x)^2}\dfrac{1}{(1+x)^2}dx
=-\sum_{k=1}^\infty\, k\, \text{Ci}(k\, \pi)
\]
and this completes the proof.

\subsection{Proof of integral 19}
In this section we treat the integral
\[
\int_0^1\dfrac{1}{\big(\pi^2+\ln(x)^2\big)^2}\dfrac{1}{(x+1)^2}
\]
We  use two methods. The first method uses differentiation of integral 12 with respect to the parameter $a$. The second method uses a a totally different method which is not very known.

\subsubsection{Differentiation to the parameter $a$}
We start with integral 12.
\[
\int_0^1\dfrac{1}{a^2+\ln(x)^2}\dfrac{1}{(1+x)^2}dx=
-\dfrac{1}{4a\pi}\psi^{(1)}\left(\dfrac{a}{2\pi}\right)+\dfrac{1}{a\pi}\psi^{(1)}\left(\dfrac{a}{\pi}\right)
\]
Differentiation with respect to the parameter $a$ and dividing both sides with the factor $-2a$ gives
\begin{multline*}
\int_0^1\dfrac{1}{\big(a^2+\ln(x)^2\big)^2}\dfrac{1}{(1+x)^2}dx=  \\
=-\dfrac{1}{8a^3\pi}\psi^{(1)}\left(\dfrac{a}{2\pi}\right)+\dfrac{1}{2a^3\pi}\psi^{(1)}\left(\dfrac{a}{\pi}\right)+
\dfrac{1}{16a^2\pi^2}\psi^{(2)}\left(\dfrac{a}{2\pi}\right)-\dfrac{1}{2a^2\pi^2}\psi^{(2)}\left(\dfrac{a}{\pi}\right)
\end{multline*}
Substitution of $a=\pi$ results in
\[
\int_0^1\dfrac{1}{\big(\pi^2+\ln(x)^2\big)^2}\dfrac{1}{(1+x)^2}dx= 
-\dfrac{1}{8\pi^4}\psi^{(1)}\left(\dfrac{1}{2}\right)+\dfrac{1}{2\pi^4}\psi^{(1)}(1)+
\dfrac{1}{16\pi^4}\psi^{(2)}\left(\dfrac{1}{2}\right)-\dfrac{1}{2\pi^4}\psi^{(2)}(1)
\]
The polygamma functions give \cite{5}
\[
\psi^{(1)}\left(\dfrac{1}{2}\right)=  3\zeta(2) \qquad \psi^{(1)}(1)=\zeta(2) \qquad
\psi^{(2)}\left(\dfrac{1}{2}\right)=-14\zeta(3)	\qquad \psi^{(2)}(1)=-2\zeta(3)
\]
Substitution gives at last
\[
\int_0^1\dfrac{1}{\big(a^2+\ln(x)^2\big)^2}\dfrac{1}{(1+x)^2}dx=\dfrac{\zeta(2)+\zeta(3)}{8\pi^4}
\]
Application of   the transformation $x \rightarrow \dfrac{1}{x}$ gives at last
\[
\int_0^1\dfrac{1}{\big(\pi^2+\ln(x)^2\big)^2}\dfrac{1}{(x+1)^2}dx=
\int_1^\infty\dfrac{1}{\big(\pi^2+\ln(x)^2\big)^2}\dfrac{1}{(x+1)^2}dx=
\dfrac{\zeta(3)+\zeta(2)}{8\pi^4}
\]
This completes the proof.

\subsubsection{Using a not so well known method}
We start with the following summation formula \cite[(1.1),(2.7)]{8}
\[
\sum_{k=1}^\infty f(k)=-2\pi\int_{-\infty}^\infty F\left(\dfrac{1}{2}+i\, t\right)
\dfrac{1}{\left(e^{\pi\, t}+e^{-\pi\, t}\right)^2}dt
\]
with $F(x)$ the primitive of the function $f(x)$. There are a lot  of conditions that the function $f(k)$ must satisfy. See\cite[section 2]{8}. The function $f(k)=\frac{1}{k^n}$ with $n>1$ certainly satisfies them.
If $f(k)=\dfrac{1}{k^3}$ the summation is known and is equal to $\zeta(3)$. The denominator in the integral can be written as a $cosh$ function. Application gives
\[
\zeta(3)=\pi\int_{-\infty}^\infty\dfrac{1}{\left(\dfrac{1}{2}+i\, t\right)^2}
\dfrac{1}{4\cosh(\pi\, t)^2}dt=\pi\int_{-\infty}^\infty\dfrac{1-4t^2-4i\, t}{(1+4t^2)^2}\dfrac{1}{\cosh(\pi\, t)^2}dt
\]
Because the left hand side of this equation is real, the right hand side should also be real. Because of that reason we omit the complex term and get
\[
\zeta(3)=\pi\int_{-\infty}^\infty\dfrac{1-4t^2}{(1+4t^2)^2}\dfrac{1}{\cosh(\pi\, t)^2}dt
\]
Splitting this equation and rewriting the $cosh$ function gives
\[
\zeta(3)=\pi\int_{-\infty}^\infty\dfrac{2}{(1+4t^2)^2}\dfrac{2}{\cosh(2\pi\, t)+1}dt-
\pi\int_{-\infty}^\infty\dfrac{1}{1+4t^2}\dfrac{2}{\cosh(2\pi\, t)+1}dt
\]
Using the transformation $u=\tanh(\pi\, t)$ gives
\[
t=\dfrac{1}{\pi}\arctanh(u) \qquad dt=\dfrac{1}{\pi}\dfrac{1}{1-u^2}du \qquad
\cosh(2\pi\, t)=\dfrac{1+\tanh(\pi\, t)^2}{1-\tanh(\pi\, t)^2}=\dfrac{1+u^2}{1-u^2}
\]
Application gives
\[
\zeta(3)=2\pi^4\int_{-1}^1\dfrac{1}{\big(\pi^2+4\arctanh(u)^2\big)^2}du-
\pi^2\int_{-1}^1\dfrac{1}{\pi^2+4\arctanh(u)^2}du 
\]
Setting
\[
\arctanh(u)=\dfrac{1}{2}\ln\left(\dfrac{1+u}{1-u}\right) \qquad x=\dfrac{1+u}{1-u} 
\qquad du=\dfrac{2}{(x+1)^2}dx
\]
gives
\[
\zeta(3)=4\pi^4\int_0^\infty\dfrac{1}{\big(\pi^2+\ln(x)^2\big)^2}\dfrac{1}{(x+1)^2}dx-
2\pi^2\int_0^\infty\dfrac{1}{\pi^2+\ln(x)^2}\dfrac{1}{(x+1)^2}dx
\]
The last integral is known from integral 13. We get
\[
\zeta(3)=4\pi^4\int_0^\infty\dfrac{1}{\big(\pi^2+\ln(x)^2\big)^2}\dfrac{1}{(x+1)^2}dx-
\dfrac{\pi^2}{6}
\]
The last term is equal to $\zeta(2)$. Rewriting the equation gives
\[
\int_0^\infty\dfrac{1}{\big(\pi^2+\ln(x)^2\big)^2}\dfrac{1}{(x+1)^2}dx=
\dfrac{\zeta(3)+\zeta(2)}{4\pi^4}
\]
This integral can be split up into two integrals when applying the transformation $x \rightarrow \dfrac{1}{x}$. 
\[
\int_0^1\dfrac{1}{\big(\pi^2+\ln(x)^2\big)^2}\dfrac{1}{(x+1)^2}dx=
\int_1^\infty\dfrac{1}{\big(\pi^2+\ln(x)^2\big)^2}\dfrac{1}{(x+1)^2}dx=
\dfrac{\zeta(3)+\zeta(2)}{8\pi^4}
\]
This completes the proof.

\section{End remarks}

We give a short overview of some of the derived integrals
\begin{align*}
&\int_0^1\dfrac{1}{\pi^2+\ln(x)^2}\dfrac{1}{(1+x)}dx=
\dfrac{1}{\pi}\sum_{k=1}^\infty \text{si}(k\pi) \\
&\int_0^1\dfrac{1}{\pi^2+\ln(x)^2}\dfrac{1}{(1+x)^2}dx=
\dfrac{1}{\pi}\sum_{k=1}^\infty k\,  \text{si}(k\pi) \\
&\int_0^1\dfrac{\ln(x)}{\pi^2+\ln(x)^2}\dfrac{1}{(1+x)}dx=
-\sum_{k=1}^\infty\text{Ci}(k\, \pi)=\dfrac{\gamma}{2}-\dfrac{1}{2}\ln(2) \\
&\int_0^1\dfrac{\ln(x)}{\pi^2+\ln(x)^2}\dfrac{1}{(1+x)^2}dx= 
-\sum_{k=1}^\infty\, k\, \text{Ci}(k\, \pi)
\end{align*}
For higher power of the second fractions of the integrals the methods of power series does not work, because in that case it is not allowed to interchange the summation and the integral.

\

There is a strong conjecture that most of the integrals can be converted in known functions.

\

Using the method of section 4.19.2 to functions like $f(k)=\dfrac{1}{k(k+1)(k+2)...}$ we can derive for example the following remarkable integral with $n \neq 0$
\[
\int_0^1\ln\left(\dfrac{(n+1)^2\pi^2+\ln(x)^2}{(n-1)^2\pi^2+\ln(x)^2}\right)\dfrac{1}{(x+1)^2}dx=
\int_1^\infty\ln\left(\dfrac{(n+1)^2\pi^2+\ln(x)^2}{(n-1)^2\pi^2+\ln(x)^2}\right)\dfrac{1}{(x+1)^2}dx=\dfrac{2}{n}
\]

This will be discussed in a future paper.

\pagebreak
\section*{Appendix: Integrals of Table 129 from Bierens de Haan.}
In this appendix we copy the integrals which are in table 129 from the book of Bierens de Haan \cite{1}. We have converted some symbols used in his table to modern symbols.

\begin{align*}
&1.\quad \int_0^1\dfrac{\ln(x)}{4\pi^2+\ln(x)^2}\dfrac{1}{1-x}dx=\dfrac{1}{4}-\dfrac{1}{2}\gamma
 \\
&2.\quad \int_0^1\dfrac{\ln(x)}{q^2+ \ln(x)^2}\dfrac{1}{1-x}dx=
\dfrac{1}{2}\left(\dfrac{\pi}{q}+\ln\left(\dfrac{2\pi}{q}\right)+\psi\left(\dfrac{q}{2\pi}\right)\right) 
\qquad\qquad\qquad\qquad\qquad\qquad
 \\
&3.\quad\int_0^1\dfrac{\ln(x)}{q^2- \ln(x)^2}\dfrac{1}{1-x}dx=
\dfrac{\pi^2}{q^2}\sum_{n=0}^\infty\dfrac{(-1)^{n-1}}{n+1}B_{2n+1}\left(\dfrac{2\pi}{q}\right)^{2n} \\
&4.\quad\int_0^1\dfrac{\ln(x)}{\left(q^2+ \ln(x)^2\right)^2}\dfrac{1}{1-x}dx=
-\dfrac{\pi^2}{q^4}\sum_{n=0}^\infty B_{2n+1}\left(\dfrac{2\pi}{q}\right)^{2n} \\
&5.\quad\int_0^1\dfrac{\ln(x)}{\left(q^2-\ln(x)^2\right)^2}\dfrac{1}{1-x}dx=
\dfrac{\pi^2}{q^2}\sum_{n=0}^\infty(-1)^{n-1}B_{2n+1}\left(\dfrac{2\pi}{q}\right)^{2n} \\
&6.\quad \int_0^1\dfrac{1}{\pi^2+ \ln(x)^2}\dfrac{1}{1+x^2}dx=\dfrac{4-\pi}{4\pi} \\
&7.\quad \int_0^1\dfrac{1}{\pi^2+ 4\ln(x)^2}\dfrac{1}{1+x^2}dx=\dfrac{1}{4\pi}\ln(2) \\
&8.\quad \int_0^1\dfrac{1}{\pi^2+16\ln(x)^2}\dfrac{1}{1+x^2}dx=
\dfrac{1}{8\pi\sqrt{2}}\left(\pi+\ln\left(\dfrac{\sqrt{2}-1}{\sqrt{2}+1}\right)\right) \\
&9.\quad \int_0^1\dfrac{1}{q^2+\ln(x)^2}\dfrac{1}{1+x^2}dx=
\dfrac{1}{4q}\left(\psi\left(\dfrac{2q+3\pi}{4\pi}\right)-\psi\left(\dfrac{2q+\pi}{4\pi}\right)\right)\\
&10.\quad \int_0^1\dfrac{\ln(x)}{\pi^2+\ln(x)^2}\dfrac{1}{1-x^2}dx=
\dfrac{1}{2}\left(\dfrac{1}{2}-\ln(2)\right) \\
&11.\quad \int_0^1\dfrac{\ln(x)}{\pi^2+4\ln(x)^2}\dfrac{1}{1-x^2}dx=\dfrac{2-\pi}{16} \\
&12.\quad \int_0^1\dfrac{\ln(x)}{\pi^2+16\ln(x)^2}\dfrac{1}{1-x^2}dx=
-\dfrac{\pi}{32\sqrt{2}}+\dfrac{1}{16}+\dfrac{1}{32\sqrt{2}}\ln\left(\dfrac{\sqrt{2}-1}{\sqrt{2}+1}\right) \\
&13.\quad \int_0^1\dfrac{\ln(x)}{\pi^2+\ln(x)^2}\dfrac{x}{1-x^2}dx=
\dfrac{1}{4}-\dfrac{1}{2}\gamma \\
&14.\quad \int_0^1\dfrac{\ln(x)}{q^2+\ln(x)^2}\dfrac{x}{1-x^2}dx=
\dfrac{1}{2}\left(\dfrac{\pi}{2q}+\ln\left(\dfrac{\pi}{q}\right)+\psi\left(\dfrac{q}{\pi}\right)\right) \\
&15.\quad \int_0^1\dfrac{\ln(x)}{q^2-\ln(x)^2}\dfrac{x}{1-x^2}dx=
\dfrac{\pi^2}{4q^2}\sum_{n=0}^\infty\dfrac{(-1)^{n-1}}{n+1}B_{2n+1}\left(\dfrac{\pi}{q}\right)^{2n} \\
&16.\quad\int_0^1\dfrac{\ln(x)}{\left(q^2+\ln(x)^2\right)^2}\dfrac{x}{1-x^2}dx=
-\dfrac{\pi^2}{4q^4}\sum_{n=0}^\infty B_{2n+1}\left(\dfrac{\pi}{q}\right)^{2n} \\
&17.\quad\int_0^1\dfrac{\ln(x)}{\left(q^2-\ln(x)^2\right)^2}\dfrac{x}{1-x^2}dx=
-\dfrac{\pi^2}{4q^4}\sum_{n=0}^\infty(-1)^{n-1} B_{2n+1}\left(\dfrac{\pi}{q}\right)^{2n}
\end{align*}


\begin{thebibliography}{1}
%
\bibitem{1}
Bierens de Haan D. {\em Nouvelles tables d'int\'egrales d\'efinies}., 1876.
%
\bibitem{3}
Blagouchin, I.V., {\em Rediscovery of Malmsten's integrals, their evaluation by contour integration methods and some related results}. The Ramanujan Journal, Volume 35, Issue 1, pp. 21-110, 2014, DOI: 10.1007/s11139-013-9528-5.
%
\bibitem{6}
Erd\'elyi, A. {\em Higher transcendental functions, Vol. I}. McGraw-Hill, 1953.  
%
\bibitem{2}
Gradshteyn, I.S., Ryzhik, I.M. {\em Tables of integrals, series and products}. Eighth edition. Academic Press 2014.
%
\bibitem{8}
Milovanovic, G.V.  {\em  Summation of series and Gaussian quadratures}, In: Approximation and Computation (R.V.M. Zahar, ed.), pp. 459-–475, Birkhäuser, 1994.

Internet: https://www.mi.sanu.ac.rs/~gvm/radovi/gvm-pur.pdf
%
\bibitem{7}
Moll,V.H. {\em Special integrals of Gradshteyn and Ryzhik, The proofs}. Vol II. CRC Press. 2016.
%
\bibitem{9} 
Saalsch\"utz, L. {\em Vorlesungen über die Bernoullischen Zahlen}. Springer, 1893.
%
\bibitem{4} 
Schl\"omilch,O. {\em Note sur quelques int\'egrales d\'efinies}.  J. reine angew. Math.,  Vol. 33,  no 316, 1846.
%
\bibitem{5}
 Weisstein, Eric W. "Polygamma Function." From MathWorld--A Wolfram Web Resource. https://mathworld.wolfram.com/PolygammaFunction.html 
\end{thebibliography}
 \end{document}